\documentclass[12pt,a4paper,twocolumn]{article}
\usepackage{amsmath,amssymb,amsfonts,amsthm}
\usepackage{mathrsfs} 
\usepackage[dvips,final]{graphicx}
\newcommand{\dfn}{\mathop{\stackrel{\mathrm{def}}{=}}}
\newcommand{\Rn}{\mathbb{R}^n}
\newcommand{\Rm}{\mathbb{R}^m}
\newcommand{\Rq}{\mathbb{R}^q}
\newcommand{\G}{\mathrm{G}}
\newcommand{\N}{\mathcal{N}}
\newcommand{\uN}{\mathrm{N}}
\newcommand{\tG}{\tilde{\mathrm{G}}}
\newcommand{\Go}{\G_{\scriptscriptstyle 0}}
\newcommand{\Gy}{G_y}
\newcommand{\Jy}{\mathcal{J}_{\scriptscriptstyle y}}
\newcommand{\Jo}{\mathcal{J}_{\scriptscriptstyle 0}}
\newcommand{\lx}{\ell\bigl(\{x_k\}\bigr)}
\newcommand{\hlx}{\widehat{\lx}}
\newcommand{\s}[1][\cdot]{\mathop{\mathrm{s}(}#1|\Gy)}
\newcommand{\elk}{\{\ell_k\}_{\scriptscriptstyle0}^{\scriptscriptstyle \uN+1}}
\newcommand{\xk}{\{x_k\}_{\scriptscriptstyle0}^{\scriptscriptstyle \uN+1}}
\newcommand{\yk}{\{y_k\}_{\scriptscriptstyle0}^{\scriptscriptstyle \uN}}
\newcommand{\hxk}{\{\hat{x}_k\}_{\scriptscriptstyle0}^{\scriptscriptstyle \uN+1}}
\newcommand{\txk}{\{\tilde{x}_k\}_{\scriptscriptstyle0}^{\scriptscriptstyle \uN+1}}
\newcommand{\hf}{\hat{f}}
\newcommand{\fk}{\{f_k\}_{\scriptscriptstyle-1}^{\scriptscriptstyle\mathrm{N}}}
\newcommand{\gk}{\{g_k\}_{\scriptscriptstyle0}^{\scriptscriptstyle \mathrm{N}}}
\newcommand{\hfk}{\{\hat{f}_k\}}

\renewcommand{\le}{\leqslant}

\theoremstyle{definition}
\newtheorem{sz_dfn}{Definition}
\theoremstyle{theorem}
\newtheorem{sz_thm}{Theorema}
\newtheorem{sz_lem}{Lemma}
\newtheorem{sz_col}{Collorary}

\begin{document}

\begin{center}
\textbf{GUARANTEED ESTIMATIONS FOR LINEAR DIFFERENCE DESCRIPTOR SYSTEMS}\\[5pt]
  Serhiy M.Zhuk\\Faculty of cybernetics\\Taras Shevchenko Kyiv
    National University, Ukraine\\beetle@unicyb.kiev.ua
\end{center}

\textbf{Abstract.}  This paper is devoted to guaranteed estimation\footnote{So-called minimax estimation} of linear functions, defined on the solutions domain of the linear descriptor difference equations (LDDE) system, where right-hand part and initial condition are arbitrary elements of the given set. Minimax estimations are build on the basis of system's state observation with unknown deterministic noise. Minimax filtration task is studied for LDDE system with special structure.
\vspace{.5cm}

\section*{Introduction.}
Suppose that vector $\xk$ satisfies linear descriptor difference equation 
\begin{equation}
  \label{eq:1}
  F_{k+1}x_{k+1}-C_k x_k = B_k f_k,k=\overline{0,\uN}
\end{equation}
with initial condition
\begin{equation}
  \label{eq:2}
  F_0x_0=B_{-1}f_{-1},
\end{equation}
where $F_k,C_k$ -- $m\times n$-matrixes, $B_k$ is $m\times p$-matrix.\\ 
We'll be interested in building minimax approximation of the linear function\footnote{$(\cdot,\cdot)_n$ denotes inner product in $\Rn$.}  $$
\lx\dfn\sum_{k=0}^{\uN+1}(\ell_k,x_k)_n,\ell_k\in\Rn,
$$ assuming that\\
\noindent\textbf{H1}\quad\emph{state $x_k$ observations are given in the form of
\begin{equation}
  \label{eq:3}
  y_k=H_k x_k+g_k,k=\overline{0,\uN},
\end{equation}
where $\gk$ is some deterministic noise},$H_k$ is $q\times n$ matrix and\\ 
\textbf{H2}\quad\emph{$\fk,\gk$  are some arbitrary elements of the ellipsoid}
\begin{equation}
  \label{eq:4}
  \begin{split}
    &\G\dfn\{(\fk,\gk):(Q_{-1}f_{-1},f_{-1})_p+\\
    &\sum_{k=0}^\uN(Q_kf_k,f_k)_p+(R_kg_k,g_k)_q\le1\},
  \end{split}
\end{equation}
where $Q_k,R_k$ are some symmetric positive-defined matrixes with appropriate dimensions.

Assume that motion of some object\footnote{Lot's of examples we can find in robototechnics~\cite{rob,karim}} is described by LDDE~\eqref{eq:1} with initial point that satisfies~\eqref{eq:2} while system disturbance ($f_k$) along with right part in~\eqref{eq:2} and noise ($g_k$) in the object's state observation model~\eqref{eq:4} are supposed\footnote{For instance, $
(\fk,\gk)$ could be measured only with some errors;  $
(\fk,\gk)$ they are random but we do not have exact information about corresponding correlation functions.} to be uncertain. Than mentioned above problem can be treated as guaranteed estimation of the object's transfer (from the set of possible initial states described by~\eqref{eq:2}), caused by uncertain disturbances, on the basis of noisy observations. Among another applications of the guaranteed estimation task studied in this paper there is a image modelling~\cite{imod} and constrained robots movement~\cite{rob}.

Let us introduce a notion of minimax a-posteriori set. At first we'll define set $\N$ as a collection of all pairs $(\xk,\fk)$ satisfying~\eqref{eq:1}-\eqref{eq:2}. Than, let us set 
\begin{equation}
  \label{Jy}
  \begin{split}
    &\Jy(\fk,\xk)\dfn \sum_{k=0}^{\uN+1}(Q_{k-1}f_{k-1},f_{k-1})_p+\\
    &(R_k(y_k-H_kx_k),y_k-H_kx_k)_q,
  \end{split}
\end{equation}
where $y_{\uN+1}=0,H_{\uN+1}=0$. If $(\fk,\gk)\in\G$ satisfies~\eqref{eq:1}-\eqref{eq:3} for some $\xk$, than $$
(\fk,\xk)\in\N, \Jy(\fk,\xk)\le1\eqno(*)
$$ and vice-versa if $(\fk,\xk)$ satisfies $(*)$ than $$
(\fk,\{g_k\dfn y_k-H_kx_k\}_0^\uN)\in\G
$$ It means that we can describe a set of all $\xk$ causing to appearance of given $\yk$ in~\eqref{eq:3} while $(\fk,\gk)$ run through some subset\footnote{This subset consists of all pairs $(\fk,\gk)\in\G$ satisfying~\eqref{eq:1}-\eqref{eq:3} for some $\xk$.} of $\G$.  
\begin{sz_dfn}
  The collection 
\begin{equation}
  \label{Gy}
  \begin{split}  
    &\Gy\dfn\{\xk|(\fk,\xk)\in\N,\\
    &(\fk,\{y_k-H_kx_k\}_0^\uN)\in\G\}
\end{split}
\end{equation}
is called \emph{a-posteriori set}. 
\end{sz_dfn}
It's obvious that real solution $\xk$ of \eqref{eq:1}-\eqref{eq:2} being observed in~\eqref{eq:3} for some $(\fk,\gk)\in\G$ belongs to $\Gy$. Hence it's naturally to look for the $\lx$ estimation \textbf{only} among the numbers from $$
L\dfn\{\lx|\xk\in\Gy\}
$$ Because of uncertain $(\fk,\gk)$ we'll use minimax strategy for finding optimal estimation $\hlx$ from within $L$: for each $\txk\in\Gy$ we need to calculate greatest distance between $\ell(\{\tilde{x}_k\})$ and $L$ -- so called \emph{guaranteed estimation error} $\sigma(\txk)$. Than we will set $\hlx=\ell(\{\tilde{x}_k\})$, where $\txk$ has a minimal $\sigma(\txk)$.
\begin{sz_dfn}
Linear function $\hlx$ is called \emph{minimax a-posteriori estimation} if 
\begin{equation*}
  \begin{split}
    &\inf_{\{\tilde{x}_k\}\in\Gy}\sup_{\{x_k\}\in\Gy}|\lx-\ell(\{\tilde{x}_k\})
    |=\\&\sup_{\{x_k\}\in\Gy}|\lx-\hlx|
  \end{split}
\end{equation*}
The non-negative number $$
\hat{\sigma}\dfn\sup_{\{x_k\}\in\Gy}|\lx-\hlx|
$$ is called \emph{minimax a-posteriori error}.
\end{sz_dfn}
In next section we shall study criteria of minimax a-posteriori estimation existence along with minimax error finiteness. It'll also be discussed a few ways of minimax estimation calculation.
\section*{Minimax a-posteriori estimation.}
\begin{sz_thm}\label{t1}
  If $\lx=\sum_{k=0}^{\uN+1}(\ell_k,x_k)_n$ and 
\begin{equation*}
\begin{split}
  &\elk\in L\dfn\{\ell_k=F'_k z_k-C'_kz_{k+1}+H'_k u_k,\\
  &F'_{\uN+1}z_{\uN+1}=\ell_{\uN+1}, z_k\in\Rm,u_k\in\Rq\}
\end{split}
\end{equation*}
than 
  \begin{align}
    &\hlx=\sum_{k=0}^{\uN+1}(\ell_k,\hat{x}_k)_n,\label{eq:hlx}\\
    &\hat{\sigma}=\bigl(1-\sum_{k=0}^\uN(y_k,R_k(y_k-H_k\hat{x}_k)_q\bigr)^
    \frac 12\bigl(\sum_{k=0}^{\uN+1}(\ell_k,p_k)\bigr)^\frac 12\label{eq:err}
  \end{align}
where $\hxk$ is a solution of 
\begin{equation}\label{eq:hxkzk}
  \begin{split}
    &F'_kz_k-C'_kz_{k+1}=H'_kR_k(y_k-H_k\hat{x}_k),k=\overline{0,\uN}\\
    &F_{k+1}\hat{x}_{k+1}-C_k\hat{x}_k=B_kQ^{-1}_kB'_kz_{k+1},\\
    &F_0\hat{x}_0=B_{-1}Q_{-1}^{-1}B'_{-1}z_0,
    F'_{\uN+1}z_{\uN+1}=0
  \end{split}
\end{equation}
and $\{p_k\}_0^{\uN+1}$ is a solution of
\begin{equation*}
  \begin{split}
    &F_{k+1}p_{k+1}=C_kp_k+B_kQ^{-1}_kB'_kd_{k+1},k=\overline{0,\uN},\\
    &F'_kd_k=C'_kd_{k+1}+\ell_k-H'_kR_kH_kp_k,\\
    &F'_{\uN+1}d_{\uN+1}=\ell_{\uN+1},F_0p_0=B_{-1}Q_{-1}^{-1}B'_{-1}d_0
  \end{split}
\end{equation*}
\end{sz_thm}
Next theorem gives a recurrence algorithm for minimax a-posteriori estimation calculation in case of special structure of matrixes $F_k,H_k$. We also suppose here that number of measurements is equal to $
\uN+1$ hence it's not necessary to set $y_{\uN+1}=0, H_{\uN+1}=0$ in~\eqref{Jy}.  
\begin{sz_thm}\label{t2}
  If $\mathop{\mathrm{rank}}
  \begin{smallmatrix}
    F_k\\H_k
  \end{smallmatrix}\equiv n$ and $B_k\equiv E$ than for any $\ell\in\Rn$ minimax a-posteriori estimation $\widehat{(\ell,x_{\uN})}$ of inner product $(\ell,x_{\uN})_n$ can be represented as $$
\widehat{(\ell,x_{\uN})}=(\ell,\hat{x}_{\uN,\uN})_n
$$ where 
\begin{equation}
    \label{eq:fltr:r}
    \begin{split}
     &\hat{x}_{k|k}=P_{k|k}F'_k(Q^{-1}_{k-1}+C_{k-1}P_{k-1|k-1}C'_{k-1})^{-1}
      \times\\
      &\times C_{k-1}\hat{x}_{k-1|k-1}+P_{k|k}H'_k R_{k}y_k,\\
      &P_{k|k}=\bigl(F'_k(Q^{-1}_{k-1}+C_{k-1}P_{k-1|k-1}C'_{k-1})^{-1}F_k+\\
      &H'_kR_kH_k\bigr)^{-1},
      P_{0|0}=(F'_0 Q F_0+H'_0 R_{0}H_0)^{-1},\\
      &\hat{x}_{0|0}=P_{0|0}H'_0 R_0y_0,
    \end{split}
  \end{equation}
\end{sz_thm}
\section*{Example.}
\begin{figure}[h]\centering
\begin{minipage}[c]{450pt} 
\includegraphics[viewport=50 200 700 700,width=450pt,clip,angle=270]{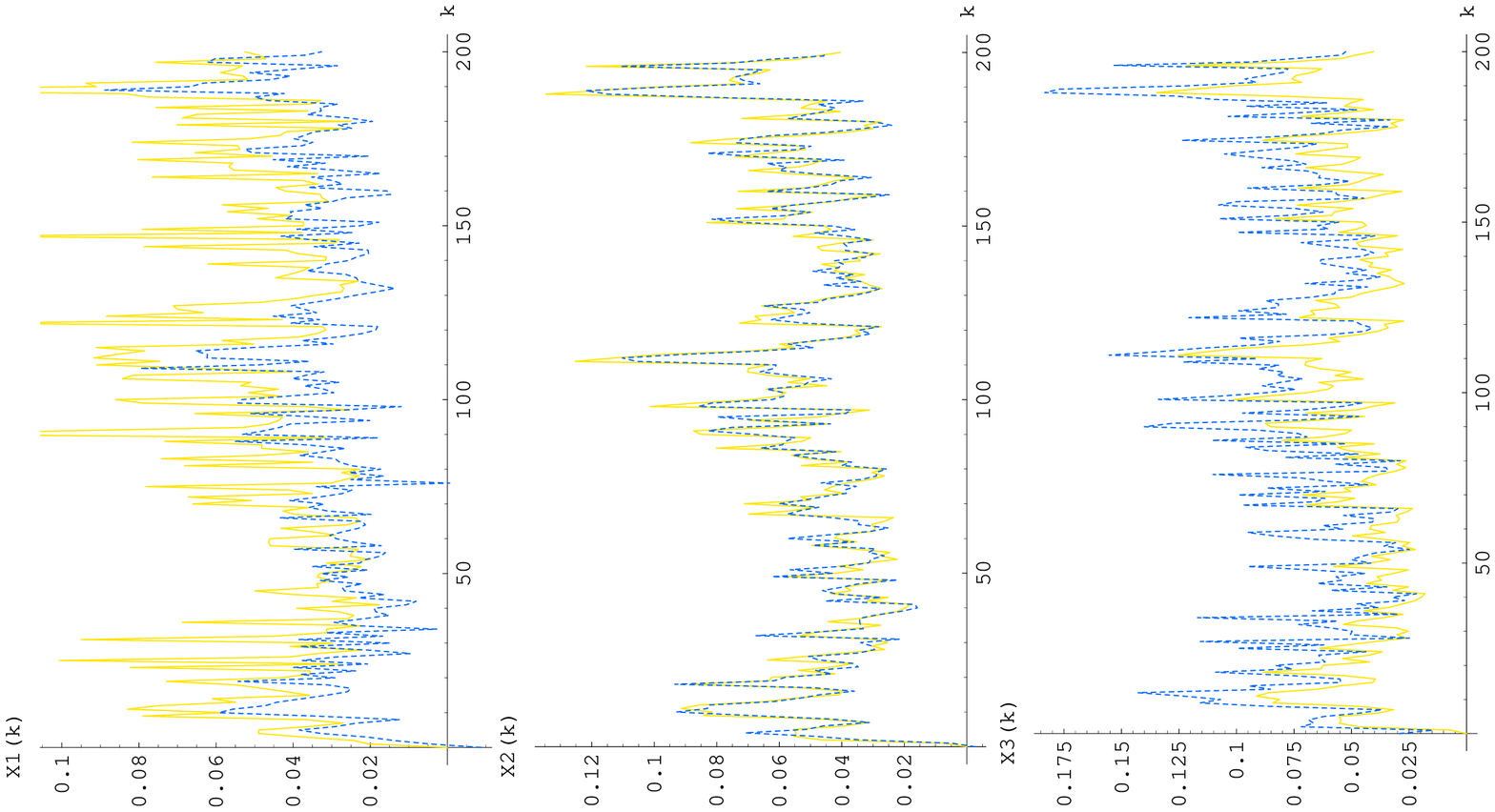}

\end{minipage}
\caption{$\hxk$(blue,dashed) and $\xk$.}
\label{fig:1}
\end{figure}

We'll illustrate theorem~\ref{t2} in case of estimating inner product $(a,x_{\uN+1})$ for linear stationary descriptor difference equation. Let us set $m=2,n=3,l=4,\uN=200$, $f_{-1}=\begin{smallmatrix}
  \frac 2{149}\\\frac 4{149}
\end{smallmatrix}
$, $F_k=
\begin{smallmatrix}
  1&&0&&0\\0&&1&&0
\end{smallmatrix}
$,$$
C_k\equiv C=
\begin{smallmatrix}
  \frac 1{40}&&0.5&&0\\0.1&&\frac 14&&\frac 3{10}
\end{smallmatrix},H_k\equiv H=
\begin{smallmatrix}
  0.6&&0.96&&2\\0.001&&2.3&&0.6\\1&&0.1&&1\\0&&0&&0.23
\end{smallmatrix}
$$ and let's choose $\fk$ and $\gk$ from the unit sphere. Simulated state $\xk$ and it's minimax a-posteriori estimation $\hxk$ are shown on the figure~\ref{fig:1}.

\section*{Proofs.}
\fontsize{8pt}{4pt}
\selectfont
\begin{proof}[Theorem's~\eqref{t1} proof.]
The following lemma gives criteria for guaranteed estimation error finiteness. 
\begin{sz_lem}
\begin{equation*}
  \begin{split}
    &\sup_{\{x_k\}\in\Gy}|\lx-\ell(\txk)|<+\infty\Leftrightarrow\elk\in L\\
  \end{split}
\end{equation*}
\end{sz_lem}
It's easy to see that $\sup_{|d|\le D}|d-c|=D+|c|$ for any real $D,c$. 
This implies to 
\begin{equation}
  \begin{split}
    &\sup_{\{x_k\}\in\Gy}|\lx-\ell(\txk)|=
    \\&\frac 12[\s[\elk]+\s[-\elk]]+
    \\&|\ell(\txk)-\frac 12[\s[\elk]-\s[-\elk]]|
  \end{split}
\end{equation}
for $\elk\in L$, so 
\begin{equation}\label{eq:hlxs}
  \begin{split}
    &\hlx=\frac 12[\s[\elk]-\s[-\elk]],\\
    &\hat{\sigma}=\frac 12[\s[\elk]+\s[-\elk]],
  \end{split}
\end{equation} 
where $\s$ is a support function of $\Gy$. Let's find $\s$.
\begin{sz_lem}
Vector $(\hfk,\hxk)$, where $\hat{f}_k=B'_kz_{k+1},k=\overline{-1,\uN}$ and $
z_k,\hat{x}_k$ are some solutions of~\eqref{eq:hxkzk}, is a minimum point of 
the $\Jy$ on $\N$
$$
\hat{J}\dfn\inf_{(\fk,\xk)\in\N}\Jy=
\sum_{k=0}^\uN(R_ky_k,y_k-H_k\hat{x}_k)_q
$$ and 
\begin{equation}
  \label{eq:JyJ0}
  \Jy(\{f_k-\hat{f}_k\}_{-1}^\uN,\{x_k-\hat{x}_k\})=
  \Jo(\fk,\xk)+\hat{J},
\end{equation}
for any\footnote{$\Jo=\Jy$ for $\yk\equiv0$.} $(\fk,\xk)\in\N$
\end{sz_lem}
If we set $ P(\fk,\xk)=\xk$ and $$
\tG\dfn\{(\fk,\xk)|\Jy(\fk,\xk)\le 1\}
$$ than 
$\Gy=P(\tG\cap\N)$ as it follows from $\Gy$ and $\G$ definitions. 
Formula~\eqref{eq:JyJ0} implies to $$\tG\cap\N=(\hfk,\hxk)+\Go\cap\N,$$ where $$
\Go\dfn\{(\fk,\xk):\Jo(\fk,\xk)\le 1-\hat{J}\}
$$ It's easy to see that $\Go=-\Go$. Hence for $\elk\in L$$$
\s[\elk]=\sum_{k=0}^{\uN+1}(\ell_k,\hat{x}_k)_n+\mathrm{s}(\elk|P(\tG_0\cap\N))
$$  Last formula with a regard to~\eqref{eq:hlxs} implies to~\eqref{eq:hlx}.

To prove~\eqref{eq:err} we need to find support function of the set $\Go\cap\N$. This task~\cite[p.164,c.16.4.1]{rkflr} is equivalent to minimisation of $\Go$-support function over affine set $ P'(\{\ell_k\})-\N^\perp$. 
It's easy to see that $$
\mathrm{Arginf s}_{P'(\{\ell_k\})-\N^\perp}\subset\mathrm{dom s}(\cdot|\Go)\cap P'(\{\ell_k\})-\N^\perp
$$ Taking into account~\cite[p.136,T.13.5]{rkflr} we can show that $$
\mathrm{dom s}(\cdot|\Go)=\{
\begin{smallmatrix}
  Q_{k-1} f_{k-1}\\H'_kR_kH_kx_k
\end{smallmatrix},k=\overline{0,\uN+1}\}\eqno(*)
$$ so regarding to structure of $L$ it's not difficult understand that minimax error can be represented as~\eqref{eq:err}.
\end{proof}
\begin{sz_col}\label{col1}
  Suppose that all assumptions of the theorem~\ref{t1} are fulfilled. To find a representation of the minimax a-posteriori estimation $\hlx$ than it is sufficient to find a solution $(\hfk,\hxk)$ of optimisation task $$
\Jy(\fk,\xk)\to\inf_{\N}
$$ 
\end{sz_col}
\begin{proof}
  As it was mentioned above we can represent minimax a-posteriori estimation in terms of a-posteriori set's support function~\eqref{eq:hlxs}.
 From the other hand each element of $\Gy$ could be treated as sum of $\hxk$ and some vector from balanced set $\N\cap\Go$. 
\end{proof}
\begin{proof}[Theorem~\eqref{t2} proof.]
  Theorem~\eqref{t2} conditions implies $L$ is equal to the whole Euclidean space $(\Rn)^{\uN+1}$, so $(0,\dots,\ell)\in L$ for any $\ell\in\Rn$ and thus we can use corollary~\ref{col1} from where and according to theorem~\eqref{t2} conditions we need to find a solution of 
  \begin{equation}
    \label{eq:opt}
    \begin{split}
      &J_{\uN+1}(\xk)\dfn\sum_{k=0}^{\uN+1}
      ||F_kx_k-C_{k-1}x_{k-1}||^2_{Q_{k-1}}+||y_k-H_kx_k||^2_{R_k},\\
      &C_{-1}=0,x_{-1}=0
    \end{split}
  \end{equation}
It is shown in~\cite{ishixara2} that we can obtain solution of~\eqref{eq:opt} 
using recurrence process~\eqref{eq:fltr:r}.
\end{proof}
\begin{proof}[Lemma's 1 proof.]
  Let $\elk\in L$. It's easy to see that 
\begin{equation*}
  \begin{split}
    &\sup_{\{x_k\}\in\Gy}|\lx-\ell(\txk)|<+\infty\Leftrightarrow\\
    &\sup_{\{x_k\}\in\Gy}|\lx|<+\infty.
  \end{split}
\end{equation*} 
On the other hand if we set $B_{-1}=B,f_{-1}=f$ than after some simple calculations we obtain
\begin{equation*}
  \begin{split}
    &\lx=\sum_{k=0}^\uN(F'_k z_k-C'_kz_{k+1},x_k)_n+(H_k x_k,u_k)_q+\\
    &(z_{\uN+1},F_{\uN+1}x_{\uN+1})_m=
    \sum_{k=0}^\uN(y_k-g_k,u_k)_q+\\
    &\sum_{k=0}^{\uN+1}(B'_{k-1}z_k,f_{k-1})_p<+\infty,\quad \forall\xk\in\Gy
  \end{split}
\end{equation*}
because of \textbf{H2}. 

Now we'll assume that $\elk\notin L$. It implies that $\lx\ne0$ for some $\xk:(\{F_kx_k,C_{k-1}x_{k-1},H_kx_k\})\equiv0$, so $\sup_{\{x_k\}\in\Gy}|\lx|=+\infty$. \end{proof}
\begin{proof}[Lemma's 2 proof.]
  Taking into account special structure of $\Jy$ and vector's projection theorem in Hilbert space~\cite{blkrshn} it's easy to prove that $\widehat{\N}\dfn\mathrm{Arginf}_{\N}\Jy\ne\varnothing$ and 
  \begin{equation*}
    \begin{split}
      \sum_{k=0}^{\uN+1}(Q_{k-1}\hf_{k-1},f_{k-1})_p+
      (H'_kR_k(y_k-H_k\hat{x}_k),x_k)_n=0,
    \end{split}
  \end{equation*}
$\forall(\fk,\xk)\in\N$ if $(\hfk,\hxk)\in\widehat{\N}$. It implies $(Q\hf,\{Q_k\hf_k\},\{H'_kR_k(y_k-H_k\hat{x}_k)\})\in\N^\perp$. On the other hand $(\hf,\hfk,\hxk)\in\N$. Now with a regard to structure of $\N$ we can show that exists $\{z_k\}$ which satisfies~\eqref{eq:hxkzk} and $\hat{f}_k=B'_kz_{k+1},k=\overline{-1,\uN}$.
\end{proof}

\end{document}